\documentclass[a4paper,12pt]{article}
\usepackage{amsmath}
\usepackage{amssymb}
\def \lim {\text {\rm lim}}

\def \tr {\text {\rm tr}}

\begin{document}
\thispagestyle{empty} \centerline{\Large{\bf Eigenvalue
inequalities for convex and log-convex functions}} \vskip 1cm
\begin{center}
{\sl Jaspal Singh Aujla\\
Department of Applied Mathematics\\
National Institute of Technology\\
Jalandhar-144011, Punjab, INDIA}\\
email: aujlajs@yahoo.com\\
and\\
{\sl Jean-Christophe Bourin\\
8 rue Henri Durel, 78510 Triel, FRANCE}\\
email: bourinjc@club-internet.fr
\end{center}

\vskip 1cm
\begin{abstract}
\noindent
We give a matrix version of the scalar inequality $f(a+b)\le f(a)+f(b)$ for positive  concave functions $f$ on $[0,\infty)$. We  show that Choi's inequality for positive unital maps and operator convex functions remains valid for monotone convex functions at the cost of unitary congruences. Some inequalities for log-convex functions
are presented and a new arithmetic-geometric mean inequality for positive matrices is given.
We also point out a simple proof of the Bhatia-Kittaneh arithmetic-geometric mean inequality.
\end{abstract}

{\bf AMS} classification: 47A30; 47B15; 15A60\\


\vspace{1cm}
\indent {\bf 1. Introduction} \\

Throughout ${\cal M}_{n}$ stands for the set of $n\times n$
complex matrices,   ${\cal H}_{n}$ for the subset of
Hermitian matrices, ${\cal S}_{n}$ for  the  positive semidefinite part of ${\cal H}_{n}$ and
${\cal P}_{n}$ for the (strictly) positive part.
We  denote by ${\cal H}_{n}(I)$ the set of $n\times n$
Hermitian matrices  with spectra in an interval $I.$

 Let $x=(x_1,x_2,\dots ,x_n)$ be an element of $\mathbb{R}^n.$
 Let $x^\downarrow$ and $x^\uparrow$ be the vectors obtained
by rearranging the coordinates of $x$ in  decreasing and
increasing order respectively. Thus $x_1^\downarrow \geq
x_2^\downarrow \geq\cdots \geq x_n^\downarrow$ and
$x_1^\uparrow\leq x_2^\uparrow \leq \cdots \leq x_n^\uparrow $.
 For  $A\in {\cal M}_n$ with real eigenvalues,
$\lambda (A)$ is a vector  of the eigenvalues
of $A$. Then, $\lambda ^\downarrow (A)$ and $\lambda ^\uparrow (A)$ can be defined
 as above.

 Let $x,y\in \mathbb{R}^n.$ The weak  majorization relation  $x\prec_wy$ means
$$
\sum_{j=1}^kx_j^\downarrow \leq \sum_{j=1}^ky_j^\downarrow,~~1\leq
k\leq n.
$$
If further equality holds for $k=n$ then we have the majorization
 $x\prec y.$ Similarly, the weak supermajorization relation
 $x\prec^wy$ means
$$
\sum_{j=1}^k x_j^\uparrow \geq \sum_{j=1}^ky_j^\uparrow,~~1\leq
k\leq n.
$$
Fan's dominance principle illustrates the relevance of
majorization in matrix theory: For $A,\, B$ in ${\cal M}_n$, the
weak majorization  $\lambda(|A|)\prec_w \lambda(|B|)$ means $\Vert
A \Vert \le \Vert B \Vert$ for all unitarily invariant norms
$\Vert\cdot\Vert$ (i.e., $\Vert UAV\Vert =\Vert A\Vert$ for all
$A$ and  all unitaries $U,V).$

Let $A=(a_{ij})$ and $B=(b_{ij})$ be elements in ${\cal M}_n.$ The  Hadamard product of $A$
and $B$ is the $n\times n$ matrix
$$
A\circ B=(a_{ij}b_{ij}).
$$
Note that $A\circ B$ appears as a principal submatrix of the Kronecker (or tensor) product
 $A\otimes B$. This is a simple but important observation of Marcus-Khan [13].
  Regarding $A\otimes B$ as an operator acting on a space ${\cal F}$
  ($\mathbb{C}^n\otimes \mathbb{C}^n$), $A\circ B$ then appears as an operator on a subspace $\cal E\subset{\cal F}$. Such a fact is expressed by saying that $A\circ B$ is the compression of $A\otimes B$ onto ${\cal E}$. A standard notation is
$
A\circ B = (A\otimes B)_{\cal E}
$.

A linear map $\Phi$ from ${\cal M}_m$ to ${\cal M}_n$ is
positive if it maps ${\cal P}_m$ to ${\cal P}_n.$ $\Phi$
 is  unital if it maps $I_m$ (identity matrix
in ${\cal M}_m$) to $I_n.$
Identifying the sets of operators on a $m$-dimensional space ${\cal F}$ and on a $n$-dimensional subspace ${\cal E}\subset{\cal F}$ with algebras ${\cal M}_m$ and ${\cal M}_n$, the compression map $\Phi(A)=A_{\cal E}$ appears as the basic example of a positive unital map. For arbitrary positive unital maps $\Phi$,
 Choi [10,11] showed that
$$
f(\Phi (A))\leq \Phi (f(A))
$$
for all operator convex functions $f$ on $I$ and all $A\in {\cal
H}_m(I).$ When $\Phi$ is a compression, this is Davis' inequality.
Choi's result is regarded as  Jensen's inequality for
noncommutative expectations.

  In Section 2 we give convexity inequalities. In particular we prove a matrix version of the basic scalar inequality
$
f(a+b)\le f(a) + f(b)
$
for positive  concave functions $f$ and $a,b\ge 0$. We also show that for monotone convex functions $f$, Choi's inequality remains valid at the cost of a unitary congruence:
 there
exists a unitary matrix $U$ such that
$$
f(\Phi (A))\leq U\Phi (f(A))U^*.
$$
When $\Phi$ is a compression, we obtain
applications to Hadamard products.

 Section 3 deals with log-convexity inequalities and  Section 4 with some related arithmetic-geometric
mean inequalities. In particular,
$$
  \prod _{j=1}^k \lambda_{j}^\uparrow(\sqrt{|AB|})\leq \prod _{j=1}^k\lambda_{j}^\uparrow
  \left (\frac{A+B}{2}\right ) ,~1\leq k \leq n
$$
for all $A,\,B\in{\cal P}_n$.  For $k=n$ this is  a classical fact: the  determinant is log-concave on ${\cal P}_n.$
\\

{\bf 2.   Eigenvalue inequalities for convex functions}\\

In [3,7] we gave several convexity inequalities for eigenvalues.
Here, we add:\\

{\bf Theorem 2.1.}
 {\it Let $f$ be a monotone concave function on $[0,\infty ) $ with $f(0)\ge 0$
  and let
 $A,B\in {\cal S}_n.$
  Then there exists unitary matrices $U$ and $V$ such that
$$
  f (A+B)\leq Uf(A)U^*+Vf(B)V^*.
$$
}
Of course, if $f$ is monotone convex with $f(0)\le 0$, we have the reverse inequality
\begin{equation}
 Uf(A)U^*+Vf(B)V^*\leq f (A+B).
\end{equation}
 Convexity (concavity) conditions are necessary
 in these results (see [1]). Theorem 2.1 yields the norm inequality (essentially Rotfel'd, see [4 p. 98])
 $$
 \Vert f(A+B)\Vert  \le \Vert f(A)\Vert  + \Vert f(B)\Vert.
 $$
 Similarly (1) implies  a trace inequality (Mc-Carthy)
 $$
 {\rm tr} A^p +{\rm tr} B^p \le {\rm tr}(A+B)^p, \quad p>1.
 $$
 Taking $f(t)=t^2, A=\left (
\begin{array}{cc}
1&0\\
0&0
\end{array}
\right )$ and $B=\left (
\begin{array}{cc}
\frac{1}{2}&\frac{1}{2}\\
\frac{1}{2}&\frac{1}{2}
\end{array}
\right )$
shows that we can not suppose $U=V$ in (1). \\

The proof of Theorem 2.1 is based on the following result [3,7]:  Let $f$ be a monotone
convex function on $I,0\in I,$  with $f(0)\leq 0.$ Then
$$
 \lambda ^\downarrow (f( XAX^*))\le \lambda ^\downarrow (X f(A)X^*)
$$
for all
 $A\in {\cal H}_n(I)$ and all contractions $X\in {\cal M}_n.$ Equivalently there exists a unitary $U$ such that
 \begin{equation}
 f( XAX^*))\le U^*X f(A)X^*U.
 \end{equation}

 {\bf Proof.}  We prove the convex version (1). We can assume that $A+B$ is invertible. Then
 $$
 A=X(A+B)X^*
 \quad {\rm and} \quad
  B=Y(A+B)Y^*
 $$
where
 $X=A^{1/2}(A+B)^{-1/2}$ and $Y=B^{1/2}(A+B)^{-1/2}$ are
 contractions. For any $T\in {\cal M}_n,$ $T^*T$ and $TT^*$ are unitarily
congruent. Hence, using (2) we have  unitaries $U_0$ and $U$ such that
 \begin{tabbing}
 \hskip 3.0cm

$f(A) $ \= $= f(X(A+B)X^*)$\\
 \>$\leq U_0 Xf(A+B)X^* U_0^*$\\
 \> $= U^*(f(A+B))^{1/2}X^*X(f(A+B))^{1/2} U,$
\end{tabbing}
so,
\begin{equation}
Uf(A)U^*\leq (f(A+B))^{1/2}X^*X(f(A+B))^{1/2}.
\end{equation}
Similarly there exists a unitary matrix $V$ such that
\begin{equation}
Vf(B)V^*\leq (f(A+B))^{1/2}Y^*Y(f(A+B))^{1/2}.
\end{equation}
Adding (3) and (4) we get
$$
 Uf(A)U^*+Vf(B)V^*\leq f (A+B)
$$
since $X^*X+Y^*Y=I_n.$ \qquad $\Box$

\vskip 10pt
 {\bf Corollary 2.2.}
 {\it Let $f$ be a   non-negative increasing concave function on $[0,\infty )$
 and let  $A,B\in {\cal S}_n.$
  Then, there exists unitary matrices $U$ and $V$ such that
$$
 Uf(A)U^*-Vf(B)V^*\leq f (|A-B|).
$$
}

 {\bf Proof.} Note that
 $$
 A\leq |A-B|+B.
 $$
 Since $f$ is increasing and concave there exists  unitaries $W,\,S,\,T$ such that
 $$
 Wf(A)W^*\leq f(|A-B|+B)\leq Sf(|A-B|)S^*+Tf(B)T^*.
 $$
Hence, we have
$$
 Uf(A)U^*-Vf(B)V^*\leq f (|A-B|)
$$
for some unitaries $U,\, V.$ \qquad $\Box$

\vskip 10pt
Other matrix versions of  basic concavity inequalities are considered in [9].
\\

Now we turn to convexity inequalities involving unital positive
maps. When necessary, ${\cal M}_n$ is identified with the algebra
$L({\cal E})$ of operators on an $n$-dimensional space ${\cal E}$.
If $X$ is an operator on a direct sum ${\cal F}=\oplus^n{\cal E}$,
then $X_{\cal E}$ stands for the compression onto the first
summand of ${\cal F}$.

Our next theorem generalizes well-known results for compressions
to arbitrary unital positive maps. The standard tool for such a
generalization is the following lemma from Stinespring's theory
[14]. The notion of a unital map on a *-subalgebra of ${\cal M}_n$
has an obvious meaning. Recall that a representation $\pi$ is a
*-homomorphism between *-subalgebras ($\pi$ preserves products,
adjoints and identities).
\\

 {\bf Lemma 2.3.} [14] {\it Let $\Phi$ be a unital  positive  map from
 a commutative *-subalgebra ${\cal A}$ of ${\cal M}_m$
to ${\cal M}_n$ identified as $L({\cal E}). $Then there exists a space ${\cal F}\supset{\cal E}$, $\dim {\cal F}\le nm$, and a
representation
$\pi$ from ${\cal A}$ to $L({\cal F})$
such that
$$
\Phi (X)=(\pi(X))_{\cal E}.
$$
}
We include a proof. It contains a simple proof of Naimark's Dilation Theorem. Say that a family of projections is total if they are mutually orthogonal and add up to the identity.

 \vskip 10pt
{\bf Proof.} ${\cal A}$ is generated by a total family of $k$
projections $E_i$, $i=1,\dots ,k$ (say $E_i$ are rank one, that is
$k=n$). Let $A_i=\Phi(E_i)$, $i=1,\dots ,n$. Since $\sum A_i$ is
the identity on ${\cal E}$, we can find operators $X_{i,j}$ such
that
$$
V=
\begin{pmatrix}
A_1^{1/2}&\dots &A_n^{1/2} \\
X_{1,1} &\dots &X_{n,1} \\
\vdots &\ddots &\vdots \\
X_{1,n-1} &\dots &X_{n,n-1}
\end{pmatrix}
$$
is a unitary operator on ${\cal F}=\oplus^n{\cal E}$. Let $R_i$ be the block matrix with the same $i$-th column than $V$ and with all other entries $0$. Then, setting $P_i=R_iR_i^*$, we obtain a total family of projections on ${\cal F}$ satifying $A_i=(P_i)_{\cal E}$. We define $\pi$ by $\pi(E_i)=P_i$.  $\Box$ \\

{\bf Theorem 2.4.} {\it Let $f$ be a  convex (resp.\ concave) function on $I$ and let $\Phi$ be a unital  positive  map from ${\cal
M}_m$ to ${\cal M}_n.$ Then
$$
 \lambda ^\downarrow (f(\Phi (A)))\prec_w \,(resp.\,\prec^w)\,\, \lambda ^\downarrow (\Phi (f(A)))
$$
for all
 $A\in {\cal H}_m(I).$ If further $f$ is also monotone then
$$
 \lambda ^\downarrow(f(\Phi (A)))\leq \,(resp.\,\ge)\,\,\lambda ^\downarrow(\Phi (f(A))).
$$
}

The proof follows from some results of [3,7] or from the following fact [8]: for $f$ convex on $I$ and  $A\in {\cal H}_n(I)$,
\begin{equation}
f(A_{\cal E}) \le \frac{U f(A)_{\cal E} U^* + V f(A)_{\cal E}
V^*}{2}
\end{equation}
for all subspaces ${\cal E}$. If further $f$ is also monotone, we can take $U=V$. \\

 {\bf Proof.} Since  $x\prec_w y$ iff $(-x)\prec^w (-y)$, it suffices to consider the convex case. Let ${\cal A}$ be the *-algebra generated by $A$. Identifying ${\cal M}_n$ with $L({\cal E})$, Lemma 2.3 yields a  representation
$\pi $ from ${\cal A}$ to $L({\cal F})$, ${\cal F}\supset{\cal E}$, such that
$$
\Phi (A)=(\pi (A))_{\cal E}.
$$
Let us denote by $\Gamma$ the compression map from $L({\cal F})$ to $L({\cal E})$. Hence, $\Phi (A)=\Gamma\circ\pi (A).$
By (5) and Fan's principle, the theorem holds for $\Gamma$. Since $\pi$ and $f$ commute, we have
$$
f(\Phi (A))=f\circ \Gamma\circ\pi (A) \prec_w \Gamma\circ  f\circ\pi (A)=\Gamma\circ \pi\circ  f (A) =\Phi (f(A))
$$
where we omitted the symbol $\lambda ^\downarrow(\cdot)$. The proof of the monotone case is similar. \qquad $\Box$

\vskip 10pt Since Hadamard products can be regarded as
compressions of tensor products and  since $|A\otimes B|=|A|
\otimes |B|$ for all $A, B$, inequality (5) for $f(t)=|t|$ gives:

\vskip 10pt
 {\bf Proposition 2.5.} {\it Let $A,B\in
{\cal H}_n $ . Then, there exist unitaries $U, V$  such that
$$
|A\circ B| \le \frac{U (|A|\circ |B|) U^* + V (|A|\circ |B|)
V^*}{2}.
$$
}

{\bf Corollary 2.6.} [12, p.\ 213] {\it For all normal matrices
$A,B\in {\cal M}_n$ and all unitarily invariant norms
$$
\Vert A\circ B\Vert \le  \Vert\, |A|\circ |B|\, \Vert.
$$
}

 The proof follows on using
$
\begin{pmatrix}
0&A^* \\ A&0
\end{pmatrix}
,
\begin{pmatrix}
0&B^* \\ B&0
\end{pmatrix}
$
in Proposition 2.5.\\

As another application of (5) we have:

\vskip 10pt
 {\bf Proposition 2.7.} {\it Let  $f$ be a submultiplicative $(f(st)\leq
f(s)f(t))$ convex function on $[ 0,\infty ).$ Then
$$
\lambda ^\downarrow (f(A\circ B))\prec_w \lambda ^\downarrow (
f(A)\circ f(B))
$$
for all $A,B\in {\cal S}_m.$ If further $f$ is also monotone then
$$
\lambda ^\downarrow (f(A\circ B))\leq \lambda ^\downarrow (
f(A)\circ f(B)).
$$
}

{\bf Proof.} Observe that the submultiplicativity of $f$ implies
$$
f(A\otimes B)\leq f(A)\otimes f(B).
$$
Let $\Gamma$ be the compression map such that $\Gamma(A\otimes B)=A\circ B$. Then, using (5) and the above inequality,
$$
f(A\circ B)= f(\Gamma (A\otimes B))
\prec_w \Gamma (f(A\otimes B))
\leq  \Gamma (f(A)\otimes f(B))
= f(A)\circ f( B)
$$
where we omitted the symbol $\lambda ^\downarrow(\cdot).$
 The monotone case can be proved similarly. \qquad $\Box$\\

We might state a version of Proposition 2.7 for supermultiplicative concave functions. The power functions are both sub and supermultiplicative and Proposition 2.7 may be applied. But the inequalities obtained follow from a stronger fact: For $A,B\in {\cal S}_m$,
$$
A^r\circ B^r \le (A\circ B)^{r}, \qquad r\in[0,1]
$$
and
$$
A^r\circ B^r \ge (A\circ B)^{r} \qquad r\in[1,2].
$$
These inequalities are special cases of Choi's inequality (see
[2]).

\vskip 10pt
We close this section by mentionning an example of positive unital map: $\Phi(A)=C\circ A$ where $C$ is a correlation matrix, i.e.\ an element of ${\cal S}_n$ with diagonal entries 1. The reader familiar with completely positive maps may note that $\Phi$ can be regarded as a compression map. \\

{\bf 3. Eigenvalue inequalities for log-convex functions }\\

Here we consider log-convexity inequalities completing [3].

\vskip 10pt
 {\bf Lemma 3.1.} {\it Let $A,B\in {\cal P}_n.$ Then
$$
\lambda ^\downarrow(\log A+\log B)\prec \lambda^\downarrow (\log
(A^{1/2}BA^{1/2}))
$$
and there exists a
unitary matrix $U$ such that
$$
\lambda_j^\downarrow (\log A+\log B) = \log \lambda_j^\downarrow
(A^{1/2}BA^{1/2}U),~~1\leq j\leq n.
$$
}

 {\bf Proof.} Let $H,K\in {\cal H}_n.$ From the Lie product formula
 $$
 e^{H+K}=\lim_{m\to \infty} (e^{H/m}e^{K/m})^{m}
 $$
  it follows (for instance [4,
Corollary IX.3.6]) that
$$
\prod _{j=1}^k\lambda_j^\downarrow (e^{H+K})\leq \prod _{j=1}^k
\lambda_j^\downarrow (e^{H/2}e^{K}e^{H/2}),~~1\leq k\leq n.
$$
The Lie formula also shows that $e^{H+K}$ and $e^Ke^H$ have the
same determinant, hence for $k=n$ equality holds. Replacing $H$
and $K$ by $\log A$ and $\log B$ we get
\begin{equation*}
\prod _{j=1}^k\lambda_j^\downarrow (e^{\log A+\log B})\leq \prod
_{j=1}^k \lambda_j^\downarrow (A^{1/2}BA^{1/2}),~~1\leq k\leq n
\end{equation*}
with equality for $k=n$. Taking logarithms  proves the first inequality. The second one follows from a famous  theorem of C.\ J.\ Thompson [15].\qquad $\Box$

\vskip 10pt
 Sometimes a more compact notation is used for inequalities involving products.
 If $x,y$ are vectors with positive coordinates the  weak log-submajorization  $x
\prec_{wlog}y$ means
$$
\prod_{j=1}^kx_j^\downarrow \leq
\prod_{j=1}^ky_j^\downarrow,~~1\leq k\leq n.
$$
Similarly, $x$ is said to be weakly log-supermajorized by $y$, in
symbol, $x\prec^{wlog}y$, if
$$
\prod_{j=1}^k x_j^\uparrow \geq \prod_{j=1}^ky_j^\uparrow,~~1\leq
k\leq n.
$$

To obtain more log-convexity inequalities, we recall:

\vskip 10pt
{\bf Theorem 3.2.} [3] {\it Let $f$ be
a convex function on $I$. Then
\begin{center}
$\lambda ^\downarrow(f(\alpha  A+(1-\alpha)B)) \prec_w
\lambda^\downarrow(\alpha f(A)+(1-\alpha )f(B))$
\end{center}
for all $A,B\in {\cal H}_{n}(I)$ and $0\leq \alpha \leq 1$.
If further $f$ is monotone, then
$$
\lambda ^\downarrow (f(\alpha A+(1-\alpha )B)\leq \lambda
^\downarrow (\alpha f(A)+(1-\alpha )f(B)).
$$
}

 {\bf Theorem 3.3.} {\it Let $f$ be a log-convex (resp.\ log-concave) function on $I$. Then, for all
 $A,B \in {\cal H}_n(I)$,
$$
\lambda ^\downarrow(f(\alpha A+(1-\alpha )B))\prec_{wlog} \,(resp.\, \prec^{wlog})\, \lambda
^ \downarrow(f(A)^{\alpha}f(B)^{1-\alpha}), \qquad 0\leq \alpha \leq 1.
$$
}

 {\bf Proof.} Since
$\log f(t)$ is   convex  on $I$, Theorem  3.2 and Lemma
3.1 yield
\begin{eqnarray*}
\lambda ^\downarrow(\log f(\alpha A+(1-\alpha )B) )& \prec_w &\lambda
^\downarrow(\alpha \log f(A)+(1-\alpha )\log f(B))\\
&=&\lambda^\downarrow(\log f(A)^{\alpha }+\log f(B)^{1-\alpha})\\
&\prec_w &\lambda ^\downarrow(\log [f(A)^{\alpha /2}f(B)^{1-\alpha
}f(A)^{\alpha /2}]).
\end{eqnarray*}
Since $\lambda ^\downarrow (f(A)^{\alpha /2}f(B)^{1-\alpha }f(A)^{
\alpha /2})=\lambda ^\downarrow (f(A)^{\alpha }f(B)^{1-\alpha })$
and $\log $ is an increasing function, we get
$$
\sum_{j=1}^k \log \lambda_j^\downarrow(f(\alpha A+(1-\alpha
)B))\leq \sum_{j=1}^k\log
\lambda_j^\downarrow(f(A)^{\alpha}f(B)^{1-\alpha}),~1\leq k\leq n
.
$$
The above inequality then implies
$$
\prod _{j=1}^k\lambda_j^\downarrow(f(\alpha A+(1-\alpha )B))\leq
\prod_{j=1}^k\lambda_j^\downarrow(f(A)^{\alpha}f(B)^{1-\alpha}),~1\leq
k\leq n.
$$
as required. If $f$ is log-concave, then
$f^{-1}(t)=\frac{1}{f(t)}$ is log-convex. Hence,  for $1\leq k\leq
n$,
$$
\prod _{j=1}^k\lambda_j^\downarrow(f^{-1}(\alpha A+(1-\alpha
)B))\leq \prod_{j=1}^k\lambda_j^\downarrow(f (A)^{-\alpha }
f(B)^{-(1-\alpha )})
$$
which implies
$$
 \prod_{j=1}^k\lambda_j^{\downarrow -1}
(f(A)^{-\alpha}f(B)^{-(1-\alpha )})\leq \prod
_{j=1}^k\lambda_j^{\downarrow -1}(f^{-1}(\alpha A+(1-\alpha )B)).
$$
Since $\lambda_j^{\downarrow -1}(H)=\lambda_j^\uparrow (H^{-1})$
and $\lambda_j^\uparrow (HK)=\lambda_j^\uparrow (KH)$ for all
$H,K\in {\cal P}_n,$ we get
$$
\prod_{j=1}^k\lambda_j^\uparrow(f(A)^{\alpha}f(B)^{1-\alpha})\leq
\prod _{j=1}^k\lambda_j^\uparrow(f(\alpha A+(1-\alpha )B)).
$$
This completes the proof. \qquad $\Box$

\vskip 3pt
For an increasing log-convex function $f$ (like $f(t)=e^t$) we can not replace in Theorem 3.3 the sign $\prec_{wlog}$ by the inequality sign (see [3]). However, we have the following statement whose proof is similar to the previous one.

\vskip 3pt
{\bf Proposition 3.4.} {\it Let $f$ be a monotone
log-convex function on $I,~0\leq \alpha \leq 1,$ and let $A,B\in
{\cal H}_n(I).$ Then there exists a unitary  matrix $U$ such that
$$
\lambda ^\downarrow (f(\alpha A+(1-\alpha )B))\leq \lambda
^\downarrow (f(A)^{\alpha /2}f(B)^{1-\alpha }f(A)^{\alpha /2}U).
$$
If $f$ is log-concave, the reverse inequality holds.
}
 \\

{\bf 4. Arithmetic-geometric mean inequalities}

\vskip 10pt
Let $A,B\in {\cal S}_n.$ Bhatia and Kittaneh [6]
proved that
\begin{equation*}
\Vert\, \sqrt{|AB|} \,\Vert \leq  \left\Vert\frac{A+B}{2}\right\Vert
\end{equation*}
for some unitarily invariant norms (for example $p$-norms for
$p\geq 2$ and trace norm) and  conjectured that it is true for all
unitarily invariant norms. Here we show the following companion
result:

\vskip 10pt

 {\bf Theorem 4.1.} {\it Let $A,B\in {\cal S}_n.$  Then, for
all $k=1, \dots,n,$
$$
\prod _{j=1}^k \lambda_j^\uparrow (\sqrt{|AB|})\leq  \prod
 _{j=1}^k \lambda_j^\uparrow \left (\frac
{A+B}{2}\right).
$$
}

Taking inverses, this theorem is equivalent to a
harmonic-geometric inequality:

\vskip 10pt
 {\bf Corollary 4.2.} {\it Let $A,\, B \in {\cal P}_n$. Then, for
all $k=1, \dots,n,$
$$
\prod_{j=1}^k \lambda_j^{\downarrow}(\sqrt{|AB|}) \ge \prod_{j=1}^k \lambda_j^{\downarrow}\left(\frac{2}{A^{-1}+B^{-1}}\right).
$$
}

We shall derive Corollary 4.2 (hence Theorem 4.1) from the
well-known result:\\

 {\bf Theorem 4.3.} (Bhatia-Kittaneh [5]) {\it Let $A,\,B\in {\cal
P}_n$. Then, there exists a unitary $U$ such that
$$
|AB| \le U\frac{A^2+B^2}{2}U^*.
$$
}

\noindent
A short proof of this theorem  follows from the following two elementary  facts in which $\Vert\cdot\Vert_{\infty}$ stands for the usual operator norm.

\vskip 10pt
{\bf Fact 1:} {\it For $A,\,B\in {\cal S}_n$ and all projections $E$,
$$
\Vert AEB\Vert_{\infty} \le \max_{\{h\in{\cal E},\ \Vert h\Vert=1\}}\Vert Ah\Vert \Vert Bh\Vert
$$
where  ${\cal E}$ stands for the range of $E$.}

\noindent Indeed, there exists rank one projection $G$ such that
$\Vert AEB\Vert_{\infty}=\Vert AEBG\Vert_{\infty}$. Letting $F$ be
the projection onto the range of $EBG$ (hence $F\le E$) we have
$$
\Vert AEB\Vert_{\infty} =\Vert AFBG\Vert_{\infty} \le \Vert AFB\Vert_{\infty}.
$$
Consequently, writing $F=h\otimes h$  we have
$$
\Vert AEB\Vert_{\infty} \le \max_{\{h\in{\cal E},\ \Vert h\Vert=1\}}\Vert Ah\Vert \Vert Bh\Vert.
$$

{\bf Fact 2:} {\it For all $A,\,B\in {\cal S}_n$ and all projections $E$ with {\rm corank}$E=k-1$,
$$
\Vert AEB\Vert_{\infty} \ge \lambda_k^\downarrow(|AB|).
$$
}
Indeed we may assume that $B$ is invertible so that there is a projection $F$, corank$F=k-1$, with $ABF=AEBF$.
Hence $\Vert AEB\Vert_{\infty}\ge \Vert ABF\Vert_{\infty}$ and Fact 2 follows from the minimax principle.

\vskip 10pt
 The proof of Theorem 4.3 is then a simple consequence of the minimax principle and the arithmetic-geometric inequality for scalars: There exists a subspace ${\cal E}$ with codimension $k-1$ such that
 \begin{align*}
 \lambda_k^\downarrow(\frac{A^2+B^2}{2})&= \max_{\{h\in{\cal E},\ \Vert h\Vert=1\}}\langle h, \frac{A^2+B^2}{2}h\rangle \\
 &\ge \max_{\{h\in{\cal E},\ \Vert h\Vert=1\}}\sqrt{\langle h, A^2h\rangle\langle h, B^2h\rangle}
 = \max_{\{h\in{\cal E},\ \Vert h\Vert=1\}}\Vert Ah\Vert \Vert
 Bh\Vert .
 \end{align*}
 Then we can apply Facts 1 and 2. \qquad $\Box$

\vskip 10pt
{\bf Proof of Corollary 4.2.} Write Bhatia-Kittaneh's inequality as
$$
(AB^2A)^{1/2} \le U\frac{A^2+B^2}{2}U^*
$$
and take inverses to get a unitary $V$ such that
$$
(A^{-1}B^{-2}A^{-1})^{1/2} \ge V\frac{2}{A^2+B^2}V^*.
$$
(Since $t\longrightarrow t^{-1}$ is not only decreasing but also operator decreasing, we can take $V=U$). Replacing $A^{-1}$ and  $B^{-1}$ by $A^{1/2}$ and $B^{1/2}$ we have
$$
(A^{1/2}BA^{1/2})^{1/2} \ge V\frac{2}{A^{-1}+B^{-1}}V^*.
$$
Since $\sqrt{|AB|}=(AB^2A)^{1/4}$, it then suffices to show that
$$
\prod_{j=1}^k\lambda_j^\downarrow((A^{1/2}BA^{1/2})^{1/2} ) \le
\prod_{j=1}^k\lambda_j^\downarrow((AB^2A)^{1/4}), \qquad k=1,
\dots ,n .
$$
But this is the same as
$$
\prod_{j=1}^k\lambda_j^{\downarrow 2}(AB) \le
\prod_{j=1}^k\lambda_j^\downarrow(A^2B^2), \qquad k=1,\dots ,n,
$$
which follows from Weyl's theorem [4, p.\ 42]. \qquad
$\Box$ \\

{\bf Proposition 4.4.} {\it Let $\{A_i\}_{i=1}^m$ be elements of
${\cal P}_n$. Then
$$
\tr (|A_1\cdots A_m|^{1/m}) \le \tr \left (\frac{A_1+\cdots
+A_m}{m}\right )
$$
 and
$$
\tr ( |A_1\cdots A_m|) \le \tr \left (\frac{A_1^m+\cdots
+A_m^m}{m}\right ).
$$
}

 {\bf Proof.} By Horn's product Theorem [12, p.\ 171]
 $$
\prod _{j=1}^{k}\lambda_{j}^\downarrow(\mid A_1\cdots A_m\mid
^{1/m})
 \le \prod _{j=1}^{k}(\lambda_{j}^\downarrow(A_1 )\cdots \lambda_{j}^\downarrow(
 A_m))^{1/m}
 $$
for all $1\leq k\leq n.$ Hence
$$
\sum _{j=1}^{k}\lambda_{j}^\downarrow(\mid A_1\cdots A_m\mid
^{1/m})
 \le \sum _{j=1}^{k}(\lambda_{j}^\downarrow(A_1 )\cdots \lambda_{j}^\downarrow(
 A_m))^{1/m}.
 $$
Then using arithmetic-geometric mean inequality for nonnegative
reals, we get
$$
\sum _{j=1}^{k}\lambda_{j}^\downarrow (\mid A_1 \cdots A_m\mid
^{1/m} )\leq \sum _{j=1}^{k}\left (\frac {
\lambda_{j}^\downarrow(A_1)+\cdots
+\lambda_{j}^\downarrow(A_m)}{m}\right ).
$$
Taking $k=n$ in the above inequality, we get the first assertion.
The second assertion follows similarly. \qquad $\Box$\\

ADDED IN PROOF:
 In a forthcoming work we will show the following
Ando-Zhan's type  inequality,
$$
||f(A+B)|| \le ||f(A)+f(B)||
$$
 for all nonnegative concave functions $f$ on $[0,\infty )$ and all $A,B  \in {\cal
 S}_n.$
The operator norm case is a striking  recent result of Tomaz
Kosem and, of course, motivates our generalization.

\end{document}